\documentclass{article}

\usepackage{amssymb,latexsym,amsmath,amsfonts,graphicx,rotating}

\newtheorem{thm}{Theorem}[section]
\newtheorem{defn}[thm]{Definition}
\newtheorem{prop}[thm]{Proposition}
\newtheorem{lem}[thm]{Lemma}
\newtheorem{rem}[thm]{Remark}

\newtheorem{ex}[thm]{Example}

\begin{document}


\begin{center}
{\bf {\Large Nearest point problem in countably normed spaces}}
\end{center}

\begin{center}
{\bf Moustafa M.  Zakaria*\footnote {$``<{\rm moustafa.m.z} \ @\ {\rm
sci.asu.edu.eg}>"$}, Nashat Faried*\footnote {$``<{\rm n}_{-}{\rm faried} \ @\ {\rm
hotmail.com}>"$} and Hany A. El-Sharkawy*\footnote  {$`` <{\rm hany.elsharkawy} \ @\
{\rm guc.edu.eg}>  <{\rm hany.elsharkawy} \ @\
{\rm sci.asu.edu.eg}>"$}}
\end{center}
\begin{center}
 {\scriptsize *
Department of Mathematics, Faculty of Science, Ain Shams University,
11566 Abbassia, Cairo, Egypt.}
\end{center}

\begin{par}
{\scriptsize
$2010$ Mathematics Subject Classification: $46A04$.

}
\end{par}

\begin{abstract}
In a countably normed space which is a linear space equipped with a 
countable number of pair-wise compatible norms, we prove the existence of a 
common nearest point (in all norms) from a point outside a nonempty subset 
if this subset is compact with respect to all norms. We also prove the uniqueness of that 
common nearest point if the completion of the space equipped with only one of its norms is uniformly convex.
\end{abstract}
\section{Introduction}
In a uniformly convex Banach space the existence and uniqueness of a nearest point 
from a point outside a given subset are guaranteed when the subset is nonempty, closed and convex.
Faried and El-Sharkawy {\rm(\cite{1})} extended this fact to the countably normed space.
They required that the completion of the space equipped with each one of all its 
norms is uniformly convex and the convex subset is closed with respect to all norms.

In this work, we prove the existence and uniqueness of the common nearest point 
from a point outside a nonempty compact subset with respect to all norms requiring only that the 
completion of the space equipped with one of its norms is uniformly convex.
 Also we give a more general theorem which guarantees the existence of a 
 nearest point if the subset is nonempty and compact with respect to all norms.
These results are useful in defining the metric projection which is a 
mapping from a given space to a given subset from the same space, 
it assigns each point in the space to its nearest point in the subset {\rm(\cite{2,3})}.
\section{Preliminaries}
\begin{defn}[Compatible norms]{\rm(\cite{4,5})}\\
{\rm Two norms $\|~\|_1$ and $\|~\|_2$ in a linear space $E$ are
said to be {\it compatible} if whenever a sequence $\{x_i\}$ in $E$
is Cauchy with respect to both norms and converges to zero with respect to one of them, it also converges to zero with respect to the other norm.}
\end{defn}

\begin{defn}[Countably normed space]{\rm(\cite{4,5})}\\
{\rm A linear space $E$ equipped with a pair-wise compatible countable family of
norms $\{\|~\|_n, n \in \mathbb{N}\}$ is said to be countably normed space. Compatibility of norms guarantees that $E$ is Hausdorff.}
If $E$ is equipped with one norm $\|~\|$, then $(E, \|~\|)$ is called {\rm a normed space}.
\end{defn}
\textbf{Remarks:}{\rm(\cite{4,5})}
\begin{itemize}
  \item Without loss of generality (by taking the equivalent system of norms $|x|_n = \max\limits_{i=1}^{n} |x|_i$ ), one can assume that the sequence of norms $\{\|~\|_n ; n = 1, 2, \ldots\}$ is increasing, i.e., $\|x\|_1 \leq \|x\|_2 \leq \cdots \leq \|x\|_n \leq \cdots , \quad \forall\ x\in E.$
  \item Any countably normed space is metrizable by using the metric $$d(x, y) = \sum\limits_{i=1}^{\infty} \frac{1}{2^i} \frac{\|x - y\|_i}{1 + \|x - y\|_i}.$$
\end{itemize}

\begin{defn}[Metric projection]{\rm(\cite{2,3})}\\
Let $K$ be a nonempty subset in a normed space $E$. The operator
$P_{K}:E\longrightarrow K$ is called a {\it metric projection
operator} if it assigns to each $x\in E$ its {\it nearest point}
$\bar{x} \in K$.
i.e. $P_{K}(x)=\bar{x}$ if $\|x-\bar{x}\|= \inf \limits_{y \in K} \|x - y\|$.
\end{defn}

\begin{defn}[Uniformly convex normed space]{\rm(\cite{6,7})}\\
{\rm A normed linear space $E$ is called {\it uniformly convex} if
for any $\varepsilon \in (0,2]$ there exists a
$\delta=\delta(\varepsilon)> 0$ such that if $x,y \in E$ with
$\|x\|=1,$ $\|y\|=1$ and $\|x-y\| \geq \varepsilon$, then
$\|\frac{1}{2}\ (x+y)\| \leq 1-\delta.$}
\end{defn}

\begin{thm}[Projection theorem in uniformly convex Banach space]{\rm(\cite{6,7})}\\
{\rm Let $K$ be a nonempty, closed and convex subset of a uniformly convex Banach space $(E, \|~\|)$. For each $x \in E\setminus K$, there exists a unique $\bar{x} \in K$ such that $$\|x -\bar{x} \|=\inf \limits_{y \in K} \|x - y\|.$$}
\end{thm}
\textbf{Notation:} For a complete countably normed space $(E, \{\|~\|_n, n \in \mathbb{N}\})$, let the completion of $E$ with respect to the norm $\|~\|_n$ be $E_n$, and for every $x \in E \subset E_n$ the norm which has been defined on $E_n$ will be $\|~\|_n$ itself.

\begin{rem}{\rm(\cite{4})}\\ Since $\|~\|_1 \leq \|~\|_2 \leq \|~\|_3 \leq \ldots$ in any countably normed space $(E, \{\|~\|_n, n \in \mathbb{N}\})$, then $E\subset \cdots \subset E_{n+1} \subset E_n \subset \cdots \subset E_1.$
\end{rem}

\begin{prop}\label{ccns}{\rm(\cite{4})}
Let $E$ be a countably normed space. Then, $E$ is complete if and only if $E = \bigcap\limits_{n = 1}^{\infty} E_n$.
\end{prop}

\begin{defn}[Uniformly convex countably normed space]\label{uccns}{\rm(\cite{1})}\\
{\rm A countably normed space $(E, \{\|~\|_n, n \in \mathbb{N}\})$ is said to be {\it uniformly
convex} if $E_n$ is uniformly convex for every $n$.}
\end{defn}

\begin{thm}\label{pcns}{\rm(\cite{1})}\\
Let $(E, \{\|~\|_n, n \in \mathbb{N}\})$ be a uniformly convex complete countably normed space, and $K$ be a nonempty convex proper subset of $E$ such that $K$ is closed in each normed space $(E_n, \|~\|_n)$. Then,
$$\forall \ x \in E\setminus K \ \ \exists ! \
\bar{x}\in K \  : \ \|x-\bar{x}\|_n = \inf \limits_{y \in K} \|x-y\|_n \quad \forall \ n \in \mathbb{N}.$$
\end{thm}

\section{Technical lemmas, remarks and propositions}
In this section, we give new lemmas and propositions which are needed to prove our main theorems.\\

Lets start with the well known proposition.
\begin{prop}
{\rm Let $(E, \|~\|)$ be a normed space. If $K \subset E$ is a nonempty compact set and
$x \in E\setminus K$, then there exists $\bar{x} \in K$ such that $\|x - \bar{x}\|=\inf \limits_{y \in K} \|x - y\|$.}

\textbf{Proof}:
Since $x\notin K$ and $K$ is
compact, then it is closed in $(E, \|~\|)$, hence
$\inf \limits_{y \in K} \|x - y\|:= d > 0 $.
Thus, there exists a sequence $\{y_i\}$ in $K$ such that
\begin{equation}
\lim \limits_{i\to \infty}
\|x - y_i\| = d > 0.
\end{equation}
It implies that $\{y_i\}$ is a bounded sequence in $K$. Then, there exists a convergent
subsequence $\{y_{i_j}\}$ to $\bar{x}$ in $K$.

$$d \le \|x -\bar{x}\| \le \|x - y_{i_j}\| +
\|y_{i_j} -\bar{x}\|$$

From equation $(1)$, $\lim \limits_{i\to \infty}
\|x - y_{i_j}\| = d$ and using $\{y_{i_j}\} \to \bar{x}$, then $d = \|x -\bar{x}\|$.
\end{prop}

\begin{rem}
$\bar{x} \in K$ is not unique in $K$. See the following example.
\end{rem}

\begin{ex}
In the space $\mathbb{R}^2$ equipped with the norm $\|(x,y)\|_{\infty}= \max \{ |x|, |y|\}$.
If $K=\{(x,y): \|(x,y)\|_{\infty} \leq 1\}$, then the set of all points in $K$ which has a minimum distance from $(2,0) \notin K$ is $K_{(2,0)}= \{(1, y): -1 \leq y \leq 1\}$.
\end{ex}
\textbf{Notation:}
 Let $K_x$ be the set of all points in $K$ which has a minimum distance from $x \notin K$.
i.e. $K_x = \{x^* \in K : \|x - x^*\|=\inf \limits_{y \in K} \|x - y\| \}$.
\\If $E$ is a linear space equipped with more than one norm, then $K_x^i$ will be the set of all points in $K$ which has a minimum distance from $x \notin K$ with respect to the norm $\|~\|_i$.

\begin{lem}
For each $x$, $K_x \subset K$ is closed in $(E, \|~\|)$, and hence it is compact.

\textbf{Proof}:
Let $\{x_i\}$ be a sequence in $K_x$ such that $\{x_i\}$ converges to $x_0 \in E$. Thus, $\|x_i - x_0\| \to 0 ~~as ~i \to \infty$.
Since $K$ is closed, then $x_0 \in K$. Moreover, $\inf \limits_{y \in K} \|x - y\| \leq \|x- x_0\| \leq \|x - x_i\| + \|x_i - x_0\| ~~ \forall i \in \mathbb{N}$.
Since $\|x - x_i\|= \inf \limits_{y\in K} \|x - y\| ~~ \forall i \in \mathbb{N}$ and $\|x_i - x_0\| \to 0 ~~as~ i \to \infty$ then $x_0 \in K_x$.
\end{lem}

At the beginning we thought that If $K$ is non empty subset of a linear space $E$ equipped with two compatible norms $\|~\|_1 \leq \|~\|_2$ and $K$ is compact with both norms, then $K^2_x \subset K^1_x$ or $K^1_x \subset K^2_x$. But this is wrong. The following example will illustrate that.

\begin{ex}
In the space $\mathbb{R}^2$ equipped with the three norms $\|~\|_{\infty} \leq \|~\|_2 \leq \|~\|_1$ where $\|(x,y)\|_{\infty}= \max \{ |x|, |y|\}, \|(x,y)\|_{2} = \sqrt{x^2 + y^2}$ and $\|(x,y)\|_{1} = |x|+ |y|$.\\
If $K=\{(x,y): \|(x,y)\|_{1} \leq 1\}$, then $K^1_{(1,1)}= \{(x, 1-x): 0 \leq x \leq 1\}$, $d_1 ((1,1),K^1_{(1,1)})= \inf \limits_{0 \leq x \leq 1} \|(1,1)-(x, 1-x)\|_1 = 1$, $K^2_{(1,1)}= \{(\frac{1}{2}, \frac{1}{2})\}$ and $d_2 ((1,1),K^2_{(1,1)})= \frac{1}{\sqrt 2}$, $K^2_{(1,1)} \subset K^1_{(1,1)}$.\\
If $K=\{(x,y): \|(x,y)\|_{\infty} \leq 1\}$, then $K^2_{(2,0)}= \{(1,0)\}$, $d_2 ((2,0),K^2_{(2,0)})=1$, $K^{\infty}_{(2,0)}= \{(1, x): -1 \leq x \leq 1\}$ and $d_{\infty} ((2,0),K^{\infty}_{(2,0)})= \max \limits_{-1 \leq x \leq 1} \|(2,0)-(1,x)\|_{\infty} = 1$, $K^2_{(2,0)} \subset K^{\infty}_{(2,0)}$.
\end{ex}

In the following two lemmas, we prove that $K^1_x \bigcap K^2_x$ is not empty by defining a new norm $\|~\|_{1+2} :=\|~\|_1 +\|~\|_2$ on $E$ and proving that $K^1_x \bigcap K^2_x = K^{1+2}_x$.

\begin{lem}
If $K$ is non a empty subset of a linear space $E$ equipped with two compatible norms $\|~\|_1 \leq \|~\|_2$ and $K$ is compact with both norms, then $K$ is non empty compact subset of the linear space $E$ equipped with the norm $\|~\|_{1+2} =\|~\|_1 +\|~\|_2$. Moreover, $(E, \|~\|_1 \leq \|~\|_2)$ is isomorphic to $(E, \|~\|_{1+2})$ by using the identity operator.

\textbf{Proof}:\\
Let $\{x_i\}$ be a bounded sequence in $(E, \|~\|_{1+2})$, then it is bounded with respect to both norms
$\|~\|_1~ and ~\|~\|_2$. Since $K$ is compact in $E$ with both norms, then there exists a convergent subsequence $\{x_{i_j}\}$ to $x_0 \in K$ with respect to $\|~\|_2$, hence it converges to the same point in $\|~\|_1$. So
$\{x_{i_j}\}$ converges to $x_0 \in K$ with respect to $\|~\|_{1+2}$.
\end{lem}

\begin{lem}
Let $E$ be a linear space equipped with two compatible norms $\|~\|_1 \leq \|~\|_2$. If $K$ is a nonempty compact subset of $E$ with both norms and $x \notin K$, then there exists a point $\bar{x} \in K$ which is a point of minimal distance from $x$ with respect to both norms.
i.e. $\|x - \bar{x}\|_1=\inf \limits_{y \in K} \|x - y\|_1$ and $\|x - \bar{x}\|_2=\inf \limits_{y \in K} \|x - y\|_2$. Moreover $K^1_x \bigcap K^2_x = K^{1+2}_x$.

\textbf{Proof}:\\
Since $K$ is a compact subset of $E$ with both norms, then it is a compact subset of $E$ with the norm $\|~\|_{1+2}= \|~\|_1 + \|~\|_2$.
Since $x \notin K$, then there exists a point $\bar{x} \in K$ such that $\|x - \bar{x}\|_{1+2}=\inf \limits_{y \in K} \|x - y\|_{1+2}$ which implies
\begin{equation}\label{2}
\|x - \bar{x}\|_1 + \|x - \bar{x}\|_2 = \inf \limits_{y \in K} \{ \|x - y\|_1 + \|x - y\|_2 \} \leq \|x - y\|_1 + \|x - y\|_2 ~~~ \forall y \in K,
\end{equation}

then

\begin{equation}\label{3}
    \|x - \bar{x}\|_1 - \|x - y\|_1 \leq  \|x - y\|_2 - \|x - \bar{x}\|_2 ~~~ \forall y \in K.
\end{equation}

Since $\|~\|_1 \leq \|~\|_2$, then

\begin{equation}\label{4}
    \|x - \bar{x}\|_1 +  \|x - y\|_1 \leq \|x - \bar{x}\|_2 + \|x - y\|_2 ~~~ \forall y \in K.
\end{equation}

Adding $(3)$ to $(4)$, we get $$\|x - \bar{x}\|_1 \leq  \|x - y\|_2 ~~~ \forall y \in K,$$ which implies
$$\|x - \bar{x}\|_1 \leq \inf \limits_{y \in K} \|x - y\|_2.$$\\
Since the last two inequalities are valid for every two norms $\|~\|_1 \leq \|~\|_2$ even if $\|~\|_2 = \|~\|_1$ or $\|~\|_2 = \frac{n+1}{n} \|~\|_1 \forall n \in \mathbb{N}$, then $\|x - \bar{x}\|_1 = \inf \limits_{y \in K}  \|x - y\|_1$. Hence $\bar{x} \in K^1_x$.\\
Now, using $(2)$, we get $$\|x - \bar{x}\|_2 = \inf \limits_{y \in K} \{ \|x - y\|_1 + \|x - y\|_2 \} - \inf \limits_{y \in K}  \|x - y\|_1 ~~~ \bar{x} \in K^{1+2}_x.$$
Hence, $\bar{x}$ must be in $K^2_x$. Thus, we  proved that $K^{1+2}_x \subset K^1_x \bigcap K^2_x$.\\

Now, we prove that $K^1_x \bigcap K^2_x \subset K^{1+2}_x$.
Assume that $\bar{x} \in K^1_x \bigcap K^2_x$, then $\|x - \bar{x}\|_{1+2} = \inf \limits_{y \in K}  \|x - y\|_1 + \inf \limits_{y \in K}  \|x - y\|_2 \leq \inf \limits_{y \in K} \|x - y\|_{1+2}$. Hence, $\bar{x} \in K^{1+2}_x$.
\end{lem}

\begin{prop}[Principle of nested sequence of compact sets]~\\
Let $(E, \{\|~\|_n, n \in \mathbb{N}\})$ be a complete countably normed space, $K_n$ be a nonempty compact subset in each normed space $(E_n, \|~\|_n )$ and $K_{n+1} \subset K_{n} ~ \forall n \in \mathbb{N}$. Then, $\bigcap \limits _{n \in \mathbb{N}} K_n$ is non empty.

\textbf{Proof}:
Since each $K_n$ is a nonempty compact subset in each normed space $(E_n, \|~\|_n ) ~ \forall n \in \mathbb{N}$ and $\|~\|_n \leq \|~\|_{n+1} ~ \forall n \in \mathbb{N}$, then $K_n$ is compact with respect all norms $\|~\|_1 \leq \ldots \leq \|~\|_n ~ \forall n \in \mathbb{N}$. Besides $\bigcap \limits_{n=m}^{\infty} K_n$ is the intersection of nonempty compact sets in the normed space $(E_m, \|~\|_m)$, then $\bigcap \limits_{n=m}^{\infty} K_n$ is nonempty in $(E_m, \|~\|_m) ~\forall m \in \mathbb{N}$.
Since $E= \bigcap \limits _{n \in \mathbb{N}} E_n$, hence $\bigcap \limits _{n \in \mathbb{N}} K_n$ is nonempty in $E$.\newpage
The following diagram illustrates that.

\begin{small}
$$ E=\bigcap \limits _{n \in \mathbb{N}} E_n \ \subset \ \ldots \ \subset \ \ \ \ \ \ \ \ E_{n+1}
\ \ \ \ \subset \ \ \ \ \ E_{n} \ \ \ \ \ \subset \ \ldots \ \subset \ E_2 \ \ \subset \ \ E_1$$
{\tiny
$$\|~\|_1\leq \|~\|_2\leq \ldots \ \ \ \ \ \ \ \ \ \ \ \ \ \ \ \|~\|_1 \leq \ldots \leq \|~\|_{n+1}\ \ \ \ \ \|~\|_1 \leq \ldots
\leq \|~\|_{n} \ \ \ \ \ \ \ \ \ \ \ \|~\|_1\leq \|~\|_2 \ \ \ \ \ \|~\|_1$$
}
$$ \ \ \ \ \ \bigcap \limits _{n \in \mathbb{N}} K_n \ \ \subset \ \ldots \ \subset \ \ \ \ \ \ \ \ K_{n+1} \ \
\ \ \subset \ \ \ \ \ K_n \ \ \ \ \subset \ \ldots \ \subset \ K_2 \ \ \subset \ \ K_1$$
\end{small}
\end{prop}

\section{Main Results}
\begin{thm}\label{mr1}
Let $(E, \{\|~\|_n, n \in \mathbb{N}\})$ be a complete countably normed space, and $K$ be a nonempty subset of $E$ such that $K$ is compact in each normed space $(E_n, \|~\|_n)$. Then,
$$\forall \ x \in E\setminus K \ \ \exists  \
\bar{x}\in K \  : \ \|x-\bar{x}\|_n = \inf \limits_{y \in K} \|x-y\|_n \quad \forall \ n \in \mathbb{N}.$$
In general $\bar{x}$ is not unique in $K$.

\textbf{Proof}:
\\Since $x\notin K$ and $K$ is
compact in each $(E_n, \|~\|_n)$, then there exists a set $K_x^n$ of points
 in $K$ which has a minimum distance with respect to $\|~\|_n$.
  Here, $K_x^n$ is the finite intersection of the sets which have a minimum distance
   with respect to the norms $\|~\|_1 \leq \|~\|_2 \leq \ldots \leq \|~\|_n$ 
   because $E_n$ is equipped with the finite norms $\|~\|_1 \leq \|~\|_2 \leq \ldots \leq \|~\|_n$.
    Using Lemma $(3.7)$, $K_x^n$ is nonempty.
\\Since $\|x\|_n \leq \|x\|_{n+1}$, $E_{n+1} \subset E_n ~~ \forall n\in \mathbb{N}$ and 
each $E_n$ has the norms $\|~\|_1 \leq \|~\|_2 \leq \ldots \leq \|~\|_n$, then $K_x^n$ are
 nested decreasing sequence of compact sets.
i.e. $K_x^{n+1} \subset K_x^n ~~ \forall n \in \mathbb{N}$. Using Proposition $(3.8)$, 
then $\bigcap \limits _{n \in \mathbb{N}} K_x^n$ is nonempty in $E$. Hence 
$\bigcap \limits _{n \in \mathbb{N}} K_x^n$ is the set of points in $K$ which have a 
minimum distance with respect to $\|~\|_n ~~ \forall n \in \mathbb{N}$.
\end{thm}

Now we use Theorem $(2.5)$ and Theorem $(4.1)$ to get the following interesting result to prove the uniqueness of the nearest point.

\begin{thm}\label{mr2}
Let $(E, \{\|~\|_n, n \in \mathbb{N}\})$ be a complete countably normed space such that for some $n_0 \in \mathbb{N}$ 
the space $(E_{n_0}, \|~\|_{n_0})$ is uniformly convex Banach space, and $K$ be a nonempty convex
 proper subset of $E$ such that $K$ is compact in each normed space 
$(E_n, \|~\|_n) \ \ \forall n \in \mathbb{N}$. Then,
$$\forall \ x \in E\setminus K \ \ \exists ! \
\bar{x}\in K \  : \ \|x-\bar{x}\|_n = \inf \limits_{y \in K}\|x-y\|_n \quad \forall \ n \in \mathbb{N}.$$

\textbf{Proof}:
\\It is the same as the previous proof with $K_x^{n_0} = \{\bar{x}\}$ because $(E_{n_0}, \|~\|_{n_0})$ is uniformly convex. The following figure illustrates that.
$$
\begin{array}{ccccccccccccc}
  E=\bigcap \limits _{n \in \mathbb{N}} E_n & \subset & \ldots & \subset & E_{n_0+1} & \subset & E_{n_0} & \subset & \ldots & \subset & E_2 & \subset & E_1 \\
  \cup & & & & \cup & & \text{\begin{rotate}{90}$\in$\end{rotate}}& & & & \cup & & \cup \\
  \{ \bar{x} \}=\bigcap \limits _{n \in \mathbb{N}} K_x^n & \subset & \ldots & \subset & K_x^{n_0+1} & \subset & \bar{x} & \subset & \ldots & \subset & K_x^2 & \subset & K_x^1
\end{array}
$$
\end{thm}

To give examples on the main results, we have to mention the following facts:

\begin{itemize}
  \item The space $\ell_p=\{ \{x_i\} : \sum \limits_{i=1}^{\infty} |x_i|^p < \infty \}$ with the
norm $\| \{x_i\} \|_p= (\sum \limits_{i=1}^{\infty} |x_i|^p)^{\frac{1}{p}}$ is uniformly convex Banach space
for each  $(1 < p < \infty)$.
  \item The space $\ell_1$ with either $\| ~ \|_1$ or $\| ~ \|_1+\| ~ \|_2$ is not uniformly convex.
  \item For each $(1 \leq p < \infty)$, the space $\ell_p$ with the norm $\| \{x_i\} \|_{\infty}=
\sup \limits_{i\in \mathbb{N}} |x_i|$ is not uniformly convex.
  \item $\ell_p \subsetneq \bigcap\limits_{q > p} \ell_{q}$.
\end{itemize}

\begin{ex}
The space $\ell_{2+0} := \bigcap\limits_{n\in \mathbb{N}} \ell_{2+\frac{1}{n}}$ is a countably normed space with 
the norms $\|~\|_3 \leq \|~\|_{2.5} \leq ... \leq \|~\|_{2+\frac{1}{n}} \leq ...$. Using proposition (\ref{ccns})
and that $\ell_{2+\frac{1}{n}}$ with the norm $ \|~\|_{2+\frac{1}{n}}$ is uniformly convex Banach space for every $n$, 
then the countably normed space $\ell_{2+0}$ is complete and uniformly convex.
$K= \{ \{x_i\} : \sum \limits_{i=1}^{\infty} |x_i|^2 \leq 1\} \subset \ell_{2+0}$ is convex and compact 
with respect to each norm. Hence $K$ satisfies theorem (\ref{pcns}).   
\end{ex}

\begin{ex}
The space $\ell_{2+0}$ with either the system of norms $\|~\|_n=\|~\|_{2+\frac{1}{n}}+\|~\|_{\infty} 
~\forall n\in \mathbb{N}$ 
or $\|~\|_{p_n}= \|~\|_{2+\frac{1}{n}}+ p_n ~\forall n\in \mathbb{N}$ where 
$p_n(\{x_i\})=\sum \limits_{i=1}^{i=n} |x_i|$ is complete countably normed space. Also $K$ is compact with 
respect to both system of norms. Hence $K$ satisfies theorem (\ref{mr1}).  
\end{ex}

\begin{ex}
The space $\ell_{2+0}$ with the system of norms $\|~\|_{p_n}= \|~\|_{2+\frac{1}{n}}+ p_n ~\forall n\in \mathbb{N}
 \setminus \{n_0\}$ and 
$\|~\|_{n_0}= \|~\|_{2+\frac{1}{n_0}}$ is complete countably normed space and ($\ell_{2+0}, \|~\|_{n_0})$ is 
uniformly convex Banach space. Hence $K$ satisfies theorem (\ref{mr2}).
\end{ex}




\begin{thebibliography}{7}

\bibitem{1} {\bf N. Faried, H. A. El-Sharkawy.},\newblock{\it~The Projection methods 
in countably normed spaces.}, Jornal of Inequalities and Application, 2015:45, (2015).

\bibitem{2} {\bf Ya. Alber},\newblock{\it~Metric and generlized projection operators in
Banach spaces}, properties and applications. In Theory and
Applications of Nonlinear Operators of Monotone and Accretive Type
(A. Kartsatos, editor). Marcel Dekker, New York 15-50, (1996).

\bibitem{3} {\bf Ya. Alber and S. Guerre-Delabriere},\newblock{\it~On the projection methods for fixed
point problems}, Analysis (Munich) 21, no. 1, 17-39, (2001).

\bibitem{4} {\bf I. M. Gel'fand, G. E. Shilov.},\newblock{\it~Generalized Functions Volume 2 Spaces of
Fundamental and Generalized Functions.}, Academic Press Inc, (1968).

\bibitem{5} {\bf  A. N. Kolmogorov and S. V. Fomin},\newblock{\it~Elements of the Theory
 of Functions and Functional Analysis},
Vol. 1 \& 2, Dover, (1999).

\bibitem{6} {\bf C. E. Chidume},\newblock{\it~Geometric Properties of Banach Spaces 
and Nonlinear Iterations}, Springer-Verlag London Limited, (2009).

\bibitem{7} {\bf B. Beauzamy},\newblock{\it~Introduction to Banach Spaces and their 
Geometry}, North-Holland Publishing Company,  (1982).

\end{thebibliography}
\end{document}